\documentclass[11pt]{amsart}
\usepackage[latin1]{inputenc}
\usepackage[cyr]{aeguill}
\usepackage[french]{babel}
\usepackage{amsxtra}
\usepackage{amssymb}
\usepackage[all]{xy}
 
\newtheorem{theoreme}{\sc Th{\'e}or{\`e}me}
\newtheorem{proposition}{\sc Proposition}
\newtheorem*{lemme}{\sc Lemme}
\newtheorem*{corollaire}{\sc Corollaire}
\theoremstyle{remark}
\newtheorem{remarque}{\it Remarque}
\newtheorem{exemple}{\it Exemple}
\newtheorem{definition}{\sc D{\'e}finition}

\font\tmsb=msbm10 at12pt
\font\smsb=msbm7
\font\ssmsb=msbm5
\newfam\msbfam
\textfont\msbfam=\tmsb
\scriptfont\msbfam=\smsb
\scriptscriptfont\msbfam=\ssmsb
\def \RM{\mathbb {R}}
\def \ZM{\mathbb{Z}}
\def \CM{\mathbb{C}}

\def \d{\partial}

\def\dt{\delta} 
\def\a{\alpha}
\def\b{\beta}
\def\e{\varepsilon}  
\def\g{\gamma}
\def\p{\varphi}

\def\lb{\left\{}
\def\rb{\right\}}
\def\G{\Gamma}   
\def \S{\Sigma}
\def \t{\tilde}

\def \to{\longrightarrow} 
\def \w{\wedge}
\def\2n{(\CM^{2n},0)}

\newcommand{\OM}{{\mathcal O}}

\newcommand{\D}{\Delta }

\medskip
\begin{document}
\title [La monodromie hamiltonienne des cycles {\'e}vanescents]
{ La monodromie hamiltonienne \\
des cycles {\'e}vanescents }
\author[Mauricio D. Garay]{ Mauricio D. Garay$^\dag$  }
\date{Mai 2005}
\address{Fachbereich Mathematik (17), Staudingerweg 9,
Johannes Gutenberg-Universit{\"a}t, 55099 Mainz, Germany}
\email{garay@mathematik.uni-mainz.de}
\thanks{$^\dag$Financ{\'e} par la bourse  GA 786/1-1 du Deutsche
Forschungsgemeinschaft}
\begin{abstract}
Nous {\'e}tudions la monodromie des cycles {\'e}vanescents d'une application holomorphe
dont les composantes sont en involution. Nous montrons sous des hypoth{\`e}ses pr{\'e}cis{\'e}es dans l'{\'e}nonc{\'e}
que le premier groupe d'homologie {\'e}vanescente non nul de la fibration
de Milnor est librement engendr{\'e} par les cycles {\'e}vanescents. 
Nous en d{\'e}duisons que l'op{\'e}rateur de variation correspondant est un
isomorphisme.
\end{abstract}
\maketitle
\parindent=0cm
{\it Cet article est d{\'e}di{\'e} {\`a} Bernard Teissier {\`a} l'occasion de
  son soixanti{\`e}me anniversaire.} 
\vskip1cm
 \section*{Introduction}
  La th{\'e}orie de Picard-Lefschetz locale {\'e}tudie les cycles {\'e}vanescents
  associ{\'e}s {\`a} un germe d'application holomorphe 
  $ f:(\CM^m,0) \to (\CM^k,0) $ et leurs monodromies (\cite{Lefschetz,Mil,Pham_monodromie}).
  Dans la plupart des cas {\'e}tudi{\'e}s, on suppose que la fibre {\`a} l'origine de $f$
  est un germe r{\'e}duit d'intersection compl{\`e}te {\`a} singularit{\'e} isol{\'e}e.\\
  Supposons que l'espace source soit muni d'une structure symplectique, i.e.,
  $\CM^m=\CM^{2n}=\{ (q,p) \}$  et consid{\'e}rons le cas o{\`u}
  les crochets de Poisson des composantes de $f$ s'annulent deux {\`a} deux:
  $$\forall i,j \leq k,\ \  \lb f_i,f_j \rb =\sum_{l=1}^n \d_{q_l} f_i \d_{p_l}f_j- \d_{p_l} f_i
  \d_{q_l}f_j=0. $$
 En g{\'e}n{\'e}ral, la fibre $ V_0 $ au dessus de l'origine d'un
  repr{\'e}sentant d'une telle application $ f $ n'est  pas
{\`a} singularit{\'e} isol{\'e}e.\\
  Prenons par exemple, $n=4,k=2$, $f_1=p_1q_1$ et
  $f_2=p_2$.
 Au voisinage de l'origine, la surface $V_0$ est analytiquement
  isomorphe au produit de la courbe plane $\G=V_0 \cap \{ q_2=0 \}$
  par un ouvert de $\CM$; le lieu singulier de $V_0$ est donc de
  dimension 1. La singularit{\'e} de $\G$ se propage dans $V_0$ le long du champ
hamiltonien  de $f_2$. Cet exemple est rigide en ce
 sens que toute d{\'e}formation de $f$ dont les composantes commutent
laisse le type analytique de $V_0$ inchang{\'e} dans un voisinage
  suffisamment petit de l'origine (\cite{lagrange}).\\

  Dans le cas o{\`u} les composantes de $f$ d{\'e}finissent un syst{\`e}me
  int{\'e}grable, i.e., pour $k=n$, on d{\'e}montre dans
  \cite{mutau} que le premier espace  de cohomologie de de
   Rham d'une fibre de Milnor de $f$ a la dimension de la base d'une
d{\'e}formation verselle Lagrangienne
  de $ V_0 $ pourvu que $ f $ admette une d{\'e}formation infinit{\'e}simalement
  Lagrangienne verselle et qu'elle v{\'e}rifie l'hypoth{\`e}se de pyramidalit{\'e} pr{\'e}cis{\'e}e
 en \ref{SS::pyramide}.\\
Ce travail est
motiv{\'e}e par cette relation entre la
g{\'e}om{\'e}trie analytique d'une part
  et la topologie d'autre part (voir {\'e}galement
\cite{Audin,Duistermaat,Nguyen}).\\

 Notons $V$ une fibre de Milnor de l'application $f$.  Sous
 l'hypoth{\`e}se de pyramidalit{\'e} qui sera pr{\'e}cis{\'e}e en \ref{SS::pyramide},
nous montrons  que
  \begin{enumerate}
\item la vari{\'e}t{\'e} $V$ est $2(n-k)$ connexe, 
\item le groupe d'homologie {\'e}vanescente $H_{2(n-k)+1}(V)$ est librement engendr{\'e}
par une base distingu{\'e}e de cycles {\'e}vanescents,
\item l'op{\'e}rateur de variation $  Var:H_{2(n-k)+1}(V,\d V) \to H_{2(n-k)+1}(V)$  est un
isomorphisme.
\end{enumerate}
Ici comme dans tout l'article, les homologies sont prises {\`a} coefficient 
dans $ \ZM $.
\section{La fibre de Milnor d'une intersection compl{\`e}te.}
\subsection{Repr{\'e}sentants standards}
Comme il est habituel en th{\'e}orie des singularit{\'e}s, nous dirons
 abusivement qu'un germe d'application  est
 une {\em intersection compl{\`e}te} si sa fibre au dessus de l'origine
 est un germe d'intersection compl{\`e}te r{\'e}duite.\\
On d{\'e}finit le {\em faisceau des applications holomorphes sur des
compacts $\CM^m$}, par limite inductive
$\OM(K)=\underrightarrow{\lim} \, \OM(U),\ K \subset U$
o{\`u} $K \subset \CM^m$ est un compact et $U \subset \CM^m$ est un ouvert.\\
Une sous vari{\'e}t{\'e} lisse {\`a} bord $M \subset \CM^m$ sera dite
{\em analytique ({\`a} bord)}
 si son  int{\'e}rieur est une vari{\'e}t{\'e} analytique complexe lisse et en
 tout point du bord
$M$ est localement  l'intersection transverse d'une vari{\'e}t{\'e} analytique
 complexe lisse avec une boule ferm{\'e}e.\\
Le bord d'une vari{\'e}t{\'e} analytique lisse compacte est une vari{\'e}t{\'e}
 analytique lisse r{\'e}elle.\\
Notons $B_\e$ la boule ferm{\'e}e de $\CM^m$,
de rayon $\e$ centr{\'e}e en l'origine. Le bord de la boule $B_\e$ est une
 sph{\`e}re que nous notons $S_\e$.
\begin{definition}
\label{D::standard}
{ Un {\em repr{\'e}sentant standard}  $\bar  f :M \to T$ d'un germe
 d'application holomorphe $  f=(f_1,\dots,f_k):(\CM^m,0) \to (\CM^k,0)$
est un repr{\'e}sentant du germe $f$ satisfaisant aux conditions suivantes:
\begin{enumerate}
\item le morphisme $\bar f$ est la restriction d'un morphisme de Stein
  $g$ d'une boule $B_\e \subset \CM^{m} $
{\`a} un voisinage compact de $g^{-1}(0)$,
\item la fibre de $ \bar f $ au dessus de $0$ est une intersection
  compl{\`e}te r{\'e}duite munie d'une stratification de
Whitney dont les strates sont transverses aux sph{\`e}res $ S_{t\e} $
pour tout $ t \in ]0,1]$,
\item il existe un syst{\`e}me fondamental de voisinages ferm{\'e}s de l'origine
  $(D_t)_{t \in ]0,1]}$ dans $T$ avec $D_1=T$,
tel qu'au dessus de $D_t$ les fibres lisses de $g$ soient transverses aux sph{\`e}res
$ S_{t\e}$.
\end{enumerate}}
\end{definition} 
\begin{remarque}{Pour tout germe $f$, il existe un repr{\'e}sentant
v{\'e}rifiant les conditions 1 et 2. En effet, l'existence d'une boule
$B_\e$ satisfaisant (2) d{\'e}coule des propri{\'e}t{\'e}s g{\'e}n{\'e}rales des
stratifications de Whitney \cite{Wh} (voir {\'e}galement \cite{Te2}).}
\end{remarque}

\subsection{Le th{\'e}or{\`e}me de section hyperplane relatif}
La proposition suivante est une variante du th{\'e}or{\`e}me de section
hyperplane de Lefschetz pour les vari{\'e}t{\'e}s analytiques compactes (au
sens d{\'e}fini dans la sous-section pr{\'e}c{\'e}dente).
\begin{proposition}
\label{P::Lefschetz}{Soient $V \subset \CM^m$ une vari{\'e}t{\'e} analytique lisse
compacte de codimension $k$,  $u:\CM^m
\to \CM$ une forme lin{\'e}aire et $P= u^{-1}(t),\ t \in \CM$ un
hyperplan. Supposons que le triplet $(V,u,t)$
satisfasse aux conditions suivantes
\begin{enumerate}
\item l'hyperplan $P$ intersecte $V$ transversalement,
\item la restriction de l'application lin{\'e}aire $u$ {\`a} $V$ est une submersion en tout point
  du bord de $V$,
\item les points critiques de la restriction de $u$ {\`a} $V$ sont des points
critiques isol{\'e}s.
\end{enumerate}
 Alors les paires $(V,V \cap P)$, $\d V,\d(V \cap P)$ sont
 respectivement $(m-k-1)$-connexe et $(m-k-2)$-connexe. En particulier, on a les annulations en homologie
$H_j(V, V \cap P)=0, H_{j-1}(\d V, \d (V \cap P))=0$ pour $0<j < m-k$.}
\end{proposition}
\begin{proof}
Il suffit de r{\'e}p{\'e}ter la d{\'e}monstration donn{\'e}e par Thom du th{\'e}or{\`e}me de
Lefschetz (\cite{Thom_ihes}).\\
Consid{\'e}rons la fonction
$ h:V \to \RM,\ \ z \mapsto \mid u(z) \mid^2  .$
Une valeur critique $ \e >0\ $ de $h$ est une valeur pour laquelle il existe
$ t \in \CM$ de module $ \sqrt \e $ tel que l'hyperplan $ u^{-1}(t) $ est tangent {\`a}
$V$. Quitte {\`a} perturber l{\'e}g{\`e}rement la fonction $h$ on peut
supposer que
\begin{enumerate}
\item les points critiques de $h$ sont des points
critiques de Morse, i.e., la forme quadratique $d^2h(z)$ est non-d{\'e}g{\'e}n{\'e}r{\'e}e
en un point critique,
\item le flot du gradient de $h$ donne une d{\'e}composition cellulaire de
la variété $V \setminus (V \cap P)$ (voir e.g. \cite{Bott,Milnor_Morse}). 
\end{enumerate}
En un point critique de $ h $, l'indice est au moins {\'e}gal -en fait on pourrait
choisir $h$ pour qu'il soit {\'e}gal- {\`a} la dimension de $ V $ .\\
Par ailleurs, en tout point du  bord de $V$, la restriction de
l'application $u$ {\`a} $V$ {\'e}tant une
submersion, le gradient de $ h $ est transverse au bord de $  V $.
Le th{\'e}or{\`e}me fondamental de la th{\'e}orie de Morse
entra{\^\i}ne que la vari{\'e}t{\'e}  $ V $ est obtenue {\`a} partir de $
V \cap P $ en y ajoutant des cellules de dimension {\'e}gale {\`a} l'indice de
chacun des points critiques de la fonction $h$
(voir e.g. \cite{Milnor_Morse}).\\ L'intersection d'une cellule de dimension
$j$ de $V$ avec le bord de $V$ est une cellule de dimension $j-1$, par cons{\'e}quent  la vari{\'e}t{\'e}  $
\d V $ est obtenue {\`a} partir de $\d(V \cap P) $ en y ajoutant des
cellules de dimension $\dim V-1$. Ceci ach{\`e}ve la d{\'e}monstration du lemme.
\end{proof}
\subsection{$N$-connexit{\'e} de la fibre de Milnor}
\label{SS::connexe}
Soit $u$ une forme lin{\'e}aire sur $\CM^m$ et $\bar f:M \to T$
un repr{\'e}sentant standard d'un germe d'intersection compl{\`e}te.\\
Consid{\'e}rons la propri{\'e}t{\'e}\\
P$(u,\bar f)$: il existe un voisinage de l'origine  $T_u \subset T$, tel que
pour toute fibre lisse $V$ de $\bar f$ au dessus de $T_u$ et tout $t \neq 0$
suffisamment proche de l'origine, le triplet $(V,u,t)$
satisfait les conditions de la Proposition \ref{P::Lefschetz}.
\begin{proposition}
\label{P::generique}
{Supposons qu'un germe $f:(\CM^m,0) \to (\CM^k,0)$  d'intersection compl{\`e}te admette
  un repr{\'e}sentant standard $\bar f:M \to T$. L'ensemble $\Omega$ des formes lin{\'e}aires sur $\CM^m$ v{\'e}rifiant la propri{\'e}t{\'e} P$(u,\bar f)$ est alors un ouvert de Zariski.}
\end{proposition}
\begin{proof}
Consid{\'e}rons l'ensemble $\Omega \subset (\CM^m)^{\spcheck}$ des formes lin{\'e}aires pour
lesquelles
\begin{enumerate}
\item la polaire relative {\`a} $f$, soit $\G$, est une vari{\'e}t{\'e} r{\'e}duite de
dimension $k$,
\item la multiplicit{\'e} d'intersection $(\G,V_0)$ {\`a} l'origine de $\G$ avec la fibre $V_0$ de $f$ au-dessus de $0$ est finie.
\end{enumerate}
L'ensemble $\Omega$ est un ouvert de Zariski (\cite{Te2}, section 4, Proposition 1).\\
Pour $u \in \Omega$, la multiplicit{\'e} d'intersection {\`a} l'origine de la
polaire $\G$ avec la fibre singuli{\`e}re $V_0$ est finie. Donc, dans un voisinage
de l'origine suffisamment petit, cette intersection est r{\'e}duite {\`a} l'origine.\\
Par cons{\'e}quent, il existe un voisinage de l'origine $T_u$ au dessus
duquel les fibres de $\bar f:M \to T$ n'intersectent pas le bord de la
polaire $\G$.\\
Au dessus de l'ouvert $T_u$, les fibres lisses du morphisme $f$
satisfont aux conditions de la Proposition \ref{P::Lefschetz}. En
effet (avec les m{\^e}mes notations),
les points critiques de la restriction de $u$ {\`a} $V$ sont les
intersections de $V$ avec $\G$.
Le nombre de ces points est donc fini.\\
Soit $t$ une valeur r{\'e}guli{\`e}re de la restriction de $u$ {\`a} $V$. L'hyperplan
$P=u^{-1}(t)$ est transverse {\`a} $V$. Par ailleurs comme le bord de $ M$
n'intersecte pas $\G$, la restriction de $u$ au bord de $V$ est une submersion.
La proposition est d{\'e}montr{\'e}e.
\end{proof}

\begin{corollaire}
\label{C::Lefschetz}
{Supposons qu'un germe  $f:(\CM^m,0) \to (\CM^k,0)$ d'intersection
  compl{\`e}te admette un repr{\'e}sentant standard $\bar f:M \to T$. Soit
  $s$ la dimension du lieu singulier de $V_0=f^{-1}(0)$.
Il existe alors, un voisinage de l'origine $T_0 \subset T$ tel que
  les fibres lisses de $\bar f$ au dessus de $T_0$ soient $(m-s-k-1)$-connexe.}
\end{corollaire}
\begin{proof}
L'intersection de $V_0=\bar f^{-1}(0) $ avec un espace vectoriel $P=\{ u=0
\}$ de codimension $s$ est une intersection compl{\`e}te $W_0=\{ (\bar f,u)=0
\}$ {\`a} singularit{\'e} isol{\'e}e. D'apr{\`e}s Milnor et Hamm(\cite{Mil,Hamm}),
la fibre de Milnor, soit $W$,
de $(f,u)$ a le type d'homotopie d'un bouquet de sph{\`e}res de dimension
$m-s-k$. Par ailleurs, la proposition pr{\'e}c{\'e}dente entra{\^\i}ne que $(V,W)$ est
$(m-s-k-1)$-connexe. Comme $W$ est $(m-s-k-1)$-connexe, le r{\'e}sultat en
d{\'e}coule.
\end{proof}
\section{Applications involutives pyramidales}
\subsection{Existence d'un repr{\'e}sentant standard dans le cas pyramidal}
\label{SS::pyramide}
Consid{\'e}rons {\`a} pr{\'e}sent le cas o{\`u} $M$ est muni d'une structure
symplectique complexe, i.e., d'une deux-forme holomorphe $\omega$
ferm{\'e}e non-d{\'e}g{\'e}n{\'e}r{\'e}e.
Le crochet de Poisson est alors d{\'e}finit par $\{ f,g \}\omega^n=df \w dg \w
\omega^{n-1}$ avec $2n=\dim M$.
\begin{definition}{Un morphisme  $f:M \to T,
\ T \subset \CM^k$ entre vari{\'e}t{\'e}s complexes est appel{\'e} {\em involutif} si les crochets de
Poisson des composantes de $f$ sont nuls.}
\end{definition}  
La d{\'e}monstration de la proposition suivante est imm{\'e}diate.
\begin{proposition}
\label{P::produit}
{ Soit $f:(\CM^{2n},0) \to (\CM^k,0)$ un germe
d'application involutive. Si la diff{\'e}rentielle {\`a} l'origine
de $f$ est de rang $j$ alors il existe un germe de symplectomorphisme $\p:(\CM^{2n},0) \to (\CM^{2n},0) $ et un
germe d'application biholomorphe $\psi:(\CM^k,0) \to (\CM^k,0)$ tels que
$$\psi \circ f \circ \p=(p_1,\dots,p_j,g)$$
 o{\`u} $g:(\CM^{2n},0) \to (\CM^{k-j},0)$ est une application involutive.}
\end{proposition}
Avec Thom (\cite{Thom_Fourier}), consid{\'e}rons l'ensemble
 $C^j(\bar f) \subset M$ des points o{\`u} une application holomorphe $\bar f:M \to T$ est de rang
 $j$.\\
L'adh{\'e}rence de $C^{k-1}(\bar f)$, $k=\dim T$, est un espace
 analytique appel{\'e} le {\em lieu critique} de $\bar f$; son image par
 $\bar f$ est appel{\'e} le {\em discriminant}. Un point du discriminant est appel{\'e}
 une {\em valeur critique}.
\begin{definition}{Une application involutive $\bar f:M \to T,\ T \subset \CM^k$ est appel{\'e}e {\em
      pyramidale} si pour tout $j \leq k$ la dimension de $C^j(\bar f)$ est au plus {\'e}gale {\`a}
      $2j$: $\dim C^j(\bar f) \leq 2j$.}
\end{definition} 
\begin{remarque}
La Proposition \ref{P::produit} montre que pour $\bar f$ est
pyramidale, l'espace $M$ est stratifi{\'e} par les strates lisses
$C^j(\bar f), j \leq k$ et que
les fibres de $\bar f$ sont stratifi{\'e}es de Whitney par les flots des
champs hamiltoniens des composantes de $\bar f$.
\end{remarque}
\begin{proposition}
\label{P::representant}
{Tout germe  $f:(\CM^{2n},0) \to (\CM^k,0)$ d'application involutive pyramidale
 admet un repr{\'e}sentant standard $\bar f:M \to T$. De plus, au dessus
  du compl{\'e}mentaire de son discriminant, l'application $\bar f$
  d{\'e}finit une fibration topologique localement triviale.}
\end{proposition}
\begin{proof}
Nous utilisons les notations de la D{\'e}finition \ref{D::standard}.\\
Soit $\bar f:M \to T$ un repr{\'e}sentant du germe $f$ satisfaisant aux conditions
$1,2$ de la D{\'e}finition \ref{D::standard}.\\
Je dis que tout point $x \in (\d V_0 \cap S_{t\e})$ admet un voisinage
$U_x \subset \CM^m$, sur lequel
les fibres lisses de l'application $f$ sont transverse {\`a} $S_{t\e}$.\\
Si $x$ est un point lisse de $V_0$ cela r{\'e}sulte du fait que la transversalit{\'e} est une condition ouverte.\\
Soit $x$ un point singulier de $V_0$. La condition de pyramidalit{\'e}
entra{\^\i}ne qu'en chaque point de la strate
de $x$, l'espace tangent est engendr{\'e} par les champs
hamiltoniens des composantes de $\bar f$; il est par cons{\'e}quent limite des
espaces tangents aux fibres lisses. En effet, tout vecteur tangent
$v_x$ en $x$ {\`a} la strate de $x$ s'{\'e}crit sous la forme
$v_x=\sum a_i X_i(x)$, $a_i \in \CM$, o{\`u} $X_i$ est le champ hamiltonien associ{\'e} {\`a} la
composante $f_i$. Le champ de vecteurs $\sum a_i X_i$ est d'une
part globalement d{\'e}fini et d'autre part tangent aux fibres de $f$;
par cons{\'e}quent, l'affirmation r{\'e}sulte {\`a} nouveau de
l'ouverture de la transversalit{\'e}.\\
Pour $t \in ]0,1]$ fix{\'e}, notons $r(t,x)$ la distance de $0$ {\`a} $\bar
    f(U_x)$. Soit $D_t$ la boule ferm{\'e}e de rayon  $r(t)=inf\{ r(t,x):x \in \d
    V_0 \}$. La restriction de $\bar f$ au dessus de $D_1$ est un
    repr{\'e}sentant standard.\\
Notons $C$ le lieu critique de $\bar f$ et $\S$ son image par $f$.\\
Soit $g:M \setminus C \to T \setminus \S$ la restriction 
de l'application $\bar f$ {\`a} $M \setminus C$. Stratifions  cette
application, en d{\'e}composant $M \setminus C$ en deux strates:
d'une part l'int{\'e}rieur de $M \setminus C$ et d'autre part le bord de $M \setminus C$.\\
Nous avons d{\'e}montr{\'e} que la condition de pyramidalit{\'e} entra{\^\i}ne
\begin{enumerate}
\item l'existence d'un
repr{\'e}sentant standard $\bar f$ pour un germe d'application pyramidale
involutive donn{\'e},
\item  que l'application $g$ v{\'e}rifie la condition
$(a_f)$ de Thom (\cite{Thom_strates}).
\end{enumerate}
Par cons{\'e}quent, le deuxi{\`e}me lemme d'isotopie de Thom appliqu{\'e} {\`a}
l'application $g:M \setminus C \to T \setminus \S$
montre que $g$ est une fibration
localement triviale  \cite{Thom_strates} (voir {\'e}galement
\cite{Gibson_Looijenga}).
Ceci conclut la d{\'e}monstration de la proposition.
\end{proof}
Nous appellerons {\em fibration de Milnor} d'un germe d'application involutive pyramidale
$f:(\CM^{2n},0) \to (\CM^k,0)$, la fibration topologique d{\'e}finie par
un repr{\'e}sentant standard de $f$ au
dessus du compl{\'e}mentaire de son discriminant.
La {\em fibre de Milnor} du germe $f$ est la vari{\'e}t{\'e}
topologique obtenue comme fibre de sa fibration de Milnor.

\subsection{$N$-connexit{\'e} de la fibre de Milnor.}
En d{\'e}pit de sa simplicit{\'e}, le th{\'e}or{\`e}me suivant va nous fournir la clef de la
d{\'e}monstration du th{\'e}or{\`e}me principal. Il marque la diff{\'e}rence
essentielle entre l'{\'e}tude des applications involutives et celle des
intersections compl{\`e}tes g{\'e}n{\'e}rales.
\begin{theoreme}
\label{T::Lefschetz} {La fibre de Milnor d'un germe d'application
involutive pyramidale $f:(\CM^{2n},0) \to (\CM^k,0)$ est une vari{\'e}t{\'e}
$2(n-k)$-connexe.}
\end{theoreme}
\begin{proof}
Soit $\bar f:M \to T$ un repr{\'e}sentant standard de $f$.
Le lieu singulier de la vari{\'e}t{\'e} $V_0=\{ \bar f=0 \}$ est l'adh{\'e}rence de
la strate $C^1(\bar f) \cap V_0$. La Proposition \ref{P::produit} entra{\^\i}ne que 
cette strate est au plus de dimension $k-1$. Donc, d'apr{\`e}s le corollaire de
la Proposition \ref{P::Lefschetz}, la fibre de Milnor de $f$ est $N$-connexe
avec $N=2n-(k-1)-k-1=2(n-k)$. Le th{\'e}or{\`e}me est d{\'e}montr{\'e}.
\end{proof}

\subsection{Points critiques quadratiques}
Un point critique d'une fonction dont le nombre de Milnor est fini
se d{\'e}compose apr{\`e}s perturbation en une somme de points critiques de
Morse ayant des valeurs critiques distinctes. Nous envisageons
l'analogue de cette propri{\'e}t{\'e} pour les application involutives.

\begin{definition}{Un point critique $ x \in M $ d'une application
involutive\\ $f:M \to T$ est appel{\'e} {\em quadratique} si il
existe un germe de symplectomorphisme $ \p:(\CM^{2n},0) \to
(\CM^{2n},x) $ et un germe d'application biholomorphe $ \psi:(\CM^k,s)  \to
(\CM^k,0) $, $ s=f(x) $, tels que
$$ \psi \circ f \circ \p(q,p)=(q_1^2+p_1^2,p_2,\dots,p_k)  .$$}
\end{definition}

Notons respectivement $ C \subset M $ et $\S \subset T$, le lieu
critique et le discriminant d'une application holomorphe $\bar f:M \to
T $.\\
D{\'e}finissons les propri{\'e}t{\'e}s suivantes: \\
(M1): l'ensemble des valeurs $ \S' \subset T$ critiques non quadratiques
forment un ensemble de codimension au moins $ 2 $ dans $ T $,\\
(M2): pour tout $ s \notin \S' $, le lieu singulier d'une fibre singuli{\`e}re $ L_s $
est connexe.\\

Nous conjecturons que ces propri{\'e}t{\'e}s sont g{\'e}n{\'e}riques,
en ce sens qu'en dehors d'un ensemble de codimension infinie (dans
l'espace vectoriel des germes d'applications involutives avec $n$ et
$k$ fix{\'e}s muni de la topologie de Whitney),
pour tout germe d'application involutive,
il existe une d{\'e}formation dont le d{\'e}ploiement est pyramidal et qui poss{\`e}de les
propri{\'e}t{\'e}s (M1) et (M2). Dans le cas $ n=1 $, une telle d{\'e}formation est
donn{\'e}e par une morsification de $ f $. Nous conjecturons {\'e}galement que la
variante r{\'e}elle de ces propri{\'e}t{\'e}s d{\'e}finit les applications involutives
stables sur une vari{\'e}t{\'e} symplectique lisse compacte.

\subsection{Cycles {\'e}vanescents d'une applications involutive
pyramidale}
\label{SS::cycles}
Soit $f:(\CM^{2n},0)  \to (\CM^k,0)$ un germe d'application involutive
pyramidale. Notons $s$ la dimension du lieu
singulier de  la fibre de $ f$ au dessus de l'origine.
Choisissons un repr{\'e}sentant standard $\bar f:M \to T$  de $f$,
des formes lin{\'e}aires $u_1,\dots,u_s:\CM^{2n} \to \CM$ de telle
sorte que
\begin{enumerate}
\item pour tout $l \leq s$, les applications $F_l=(\bar
  f,u_1,\dots,u_{l-1}) $ et $F_1=f$ v{\'e}rifient la propri{\'e}t{\'e}
  P$(u_l,F_l)$ (voir section \ref{SS::connexe}),
\item la vari{\'e}t{\'e}  $V_0 \cap P_0$ est une intersection
compl{\`e}te {\`a} singularit{\'e} isol{\'e}e dont l'id{\'e}al est engendr{\'e} par les
composantes de $f$ et par les $u_i$.
\end{enumerate}
L'existence de $\bar f$ et de $u$ v{\'e}rifiant la premi{\`e}re condition se
d{\'e}montre comme dans la Proposition 6.
La deuxi{\`e}me condition est v{\'e}rifi{\'e}e g{\'e}n{\'e}riquement; il suffit en effet
que $V_0 \cap P_0$ soit un espace analytique r{\'e}duit et que
l'intersection de $P_0$ avec le lieu singulier de $V_0$ soit {\'e}gal
{\`a} l'origine (avec {\'e}ventuellement une certaine multiplicit{\'e}).\\
Supposons {\`a} pr{\'e}sent
que l'application $\bar f$ v{\'e}rifie les conditions (M1) et (M2).\\
La fibre au dessus de l'origine de l'application $(f,u):(\CM^{2n},0) \to
(\CM^{k+s},0)$, $u=(u_1,\dots,u_s)$, est un germe d'intersection compl{\`e}te {\`a}
singularit{\'e} isol{\'e}e. La fibre de Milnor $V \cap P$ de $(f,u)$ est l'intersection de la fibre de Milnor $V$ de $f$ avec un espace affine $P$.\\
On sait construire un syst{\`e}me distingu{\'e} de cycles {\'e}vanescents
$\g_1,\dots,\g_\mu$ qui engendrent le groupe $H_{2(n-k)+1}(V \cap P) \simeq \ZM^\b$,
$\b \leq \mu$ (\cite{Hamm,Mil}).\\
Ces cycles {\'e}vanescents sont  associ{\'e}s aux intersections du
discriminant de $(\bar f,u)$ avec une droite g{\'e}n{\'e}rique $L \subset
\CM^{k+s}$.\\
L'inclusion $V\cap P \subset V$ donne lieu {\`a} une application
naturelle en homologie
$$\xi:H_{2(n-k)+1}(V \cap P) \to H_{2(n-k)+1}(V).$$
On peut diviser les cycles {\'e}vanescents de $(f,u)$ en deux groupes suivant qu'ils
correspondent {\`a} une valeur singuli{\`e}re ou bien r{\'e}guli{\`e}re de $f$.
L'image par $\xi$ d'un cycle {\'e}vanescent correspondant {\`a}
une valeur r{\'e}guli{\`e}re de $f$ est {\'e}videmment nulle.\\
A priori, il est possible que plusieurs cycles {\'e}vanescents de $V \cap P$ s'{\'e}vanouissent pour la m{\^e}me
valeur critique de $f$. On choisit arbitrairement l'un d'entre eux.
On note $  \g_1,\dots, \g_l $ les cycles {\'e}vanescents de $V \cap P$
ainsi obtenus.\\
Le nombre $l$ et $\mu$ sont respectivement {\'e}gaux {\`a} la multiplicit{\'e} du discriminant de
$f$ et de $(f,u)$ {\`a} l'origine.
\begin{definition}{On appelle {\em base distingu{\'e}e de cycles
      {\'e}vanescents} de la fibre de Milnor $V$
d'un germe $f:(\CM^{2n},0) \to (\CM^k,0)$ d'application pyramidale involutive
(associ{\'e} au choix d'une droite $L \subset T$ et d'une base distingu{\'e}e de
chemins dans cette droite), les images des cycles $  \g_1,\dots, \g_l
$ pr{\'e}c{\'e}demment construits par
l'application $\xi:H_{2(n-k)+1}(V \cap P) \to H_{2(n-k)+1}(V) $.}
\end{definition}

\subsection{{\'E}nonc{\'e} du th{\'e}or{\`e}me principal}

\begin{theoreme}
\label{T::principal}
{Soit $f:(\CM^{2n},0)  \to (\CM^k,0)$ un germe
d'application pyramidale involutive,
v{\'e}rifiant les conditions $(M1)$ et $(M2)$. Le
groupe d'homologie $ H_{2(n-k)+1}(V) $ de la fibre de Milnor $V$ de $f$ est librement engendr{\'e} par toute base distingu{\'e}e de cycles {\'e}vanescents.
En particulier, le rang de ce groupe est {\'e}gal {\`a} la multiplicit{\'e} {\`a} l'origine du discriminant de $f$.}
  \end{theoreme}
\begin{exemple}{D{\'e}finissons l'application pyramidale $f:(\CM^6,0) \to (\CM^2,0)$ par
 $f_1=p_1q_1+p_2q_2$, $f_2=p_2q_2+p_3q_3$. Le discriminant de $f$ est
    le germe {\`a} l'origine de la r{\'e}union de trois droites $\{ s_1=0 \} \cup \{s_2=0 \} \cup \{
s_1=s_2 \}$. On a par cons{\'e}quent $H_1(V)=H_2(V)=0$ et $H_3(V) = \ZM^3$.}
\end{exemple}
\begin{exemple}
On peut g{\'e}n{\'e}raliser l'exemple pr{\'e}c{\'e}dent de la fa{\c c}on suivante. Soit
$g=(g_1,\dots,g_n):(\CM^{2n},0) \to (\CM^n,0)$ le germe d'application d{\'e}finie
par $g_i=p_i q_i$ et $R:\CM^n \to \CM^k$ une application lin{\'e}aire
surjective. Consid{\'e}rons l'application involutive $f=R g:(\CM^{2n},0)
\to (\CM^k,0)$. La
fibre de Milnor de $f$ est $2(n-k)$ connexe et son nombre de Betti
$b_{2(n-k)+1}$ est {\'e}gal au coefficient binomial
$\begin{pmatrix} n  \\ k-1 \end{pmatrix}$.
\end{exemple}
\begin{corollaire}{Sous les hypoth{\`e}ses du Th{\'e}or{\`e}me \ref{T::principal},
  l'op{\'e}rateur de variation $Var:H_{2(n-k)+1}(V,\d V) \to H_{2(n-k)+1}(V)$
est un isomorphisme.}
\end{corollaire}
\begin{proof}
Commen{\c c}ons par montrer que l'op{\'e}rateur de variation est bien d{\'e}fini
en dépit du fait que la fibre $V_0$ de $ f $ au-dessus de $ 0 $ ne soit pas {\`a} 
  singularit{\'e} isol{\'e}e.\\ 
La Proposition \ref{P::Lefschetz} et les suites exactes des
  couples $(V,\d V)$ et $(V \cap P, \d(V \cap P))$ donnent lieu {\`a} un
  diagramme commutatif (ici et dans la suite de l'article, nous
   utilisons les notations habituelles pour les morphismes injectifs, surjectifs et les isomorphismes),
$$ \xymatrix{  
H_{j+1}(V) \ar[r]& H_{j+1}(V,\d V) \ar[r]& H_{j}(\d V) \ar[r]& H_{j}(V) \\
H_{j+1}(V \cap P) \ar[r] \ar@{->>}[u]&  H_{j+1}(V
              \cap P,\d (V \cap P)) \ar[r] \ar[u] & H_{j}(\d (V \cap P))
              \ar@{->>}[u] \ar[r]&  H_{j}(V \cap P)
              \ar[u]_\simeq} $$
avec $j=2(n-k)$.\\
Donc toute classe d'homologie relative dans $H_{2(n-k)+1}(V , \d V)$
peut {\^e}tre repr{\'e}sent{\'e}e par un {\'e}l{\'e}ment de $H_{2(n-k)+1}(V \cap P , \d (V \cap P))$.\\ 
Soit $ \dt \subset V \cap P $ un cycle dont le bord est dans $\d(V
\cap P)$. Ce cycle d{\'e}finit un {\'e}l{\'e}ment d'homologie relative $ [\dt]
\in H_{2(n-k)+1}(V,\d V) $.
Notons $  \sigma $ un chemin qui parcourt le bord de $ S $ dans le sens
  direct.\\
La vari{\'e}t{\'e} $V_0 \cap P$ est {\`a} singularit{\'e} isol{\'e}e par cons{\'e}quent le
  chemin $ \sigma $ se rel{\`e}ve en un diff{\'e}omorphisme $ h $ de $ V \cap P $ d{\'e}finit {\`a}
  isotopie pr{\`e}s qui fixe le bord de $V \cap P$ (voir
  e.g. \cite{Looijenga}). On définit l'opérateur de variation par
$$Var: H_{2(n-k)+1}(V,\d V) \to  H_{2(n-k)+1}(V),\ \ [\dt]
  \mapsto  [h(\dt)-\dt].$$
Prenons pour base de  $H_{2(n-k)+1}(V,\d V)$ la base duale d'une base
distingu{\'e}e de  $H_{2(n-k)+1}(V)$.
Dans ces bases, la formule de Picard-Lefschetz entra{\^\i}ne que la
  matrice de l'op{\'e}rateur de variation est une matrice triangulaire
  dont les éléments diagonaux sont égaux à 1 ou bien à -1.
(voir e.g. \cite{AVGII}, Chapitre I, section 2.5).
Ce qui d{\'e}montre le corollaire.  
\end{proof}

\section{D{\'e}monstration du Th{\'e}or{\`e}me principal.}
Nous utilisons les notations de \ref{SS::cycles}.
\begin{lemme}{Soit $ \a,\b \in H_{2(n-k)+1}(V) $ les images de deux
 cycles de $V \cap P$ qui sont
{\'e}vanescents en une m{\^e}me valeur critique $s \in L \cap \S$ alors $ \a=\pm \b $.}
\end{lemme}
\begin{proof}
Notons $ \D $ le lieu singulier de la fibre singuli{\`e}re $ V_{z}=\bar f^{-1}(z) $, $ z \in L
\cap \S $.
Le lemme est une propri{\'e}t{\'e} locale  au voisinage de $ \D $,
on peut donc se restreindre {\`a} un voisinage tubulaire de $ \D $ dans $ M $.
\\
Les hypoth{\`e}ses (M1) et (M2) entra{\^\i}nent que
la vari{\'e}t{\'e} $ \D $ est lisse, connexe de dimension $k-1$.
De plus, l'hypoth{\`e}se de
pyramidalit{\'e} entra{\^\i}ne que l'espace tangent {\`a} $\Delta$ en chaque point
est engendr{\'e} par  les champs hamiltonien des composantes de $ \bar f$.\\
Les cycles $ \a,\b $ sont associ{\'e}s {\`a} des points critiques $ x,x' \in
V_z $ de m{\^e}me valeur critique $z$, i.e., $ \bar f(x)=\bar f(x')=z $.
Quitte {\`a} renum{\'e}roter les $ \bar f_i $ on peut supposer que les
hamiltoniens $X_1,\dots,X_{k-1} $ des fonctions
$ \bar f_1,\dots \bar f_{k-1} $ sont lin{\'e}airement ind{\'e}pendants en chaque point
d'un voisinage $ U \subset M $ de $ x \in \D$.\\
Commen{\c c}ons par d{\'e}montrer le lemme dans le cas o{\`u} $ x' \in U $.\\
Notons $ \p_i^t:\Omega \to U $ le flot du champ
hamiltonien de $ \bar f_i $ au temps $ t \in \CM $ et
consid{\'e}rons l'application
$$ \Phi:\Omega \to U,\ \ (t_1,\dots,t_{k-1}) \mapsto
(\p_1^{t_1}(x),\dots,\p_{k-1}^{t_{k-1}}(x)), $$
le domaine $ \Omega \subset \CM^{k-1}$ {\'e}tant choisi de telle sorte que
l'application $\Phi$ soit surjective.\\
Choisissons un chemin analytique dans $ \Omega$ dont l'image par $ \Phi $ est un chemin plong{\'e}
$C_{x,x'}$ joignant $ x $ {\`a} $ x'$.
{\`a} l'aide des champs hamiltoniens $ X_1,\dots,X_{k-1} $ construisons
un champ hamiltonien $X$ tangent {\`a} $ C_{x,x'} $ d{\'e}finit dans un
voisinage tubulaire du chemin $ C_{x,x'} \subset M$.\\
Le flot du champ $ X $ donne une homotopie entre des repr{\'e}sentants des cycles $ \a $ et $\pm \b
$ (l'ambigu{\"\i}t{\'e} du signe provient du fait que le cycle {\'e}vanescent est d{\'e}fini de fa{\c c}on unique {\`a} l'orientation pr{\`e}s).
Le lemme est d{\'e}montr{\'e} pour $ x' $ voisin de $ x $.\\
Dans le cas g{\'e}n{\'e}ral, on d{\'e}coupe le chemin joignant $ x $ {\`a} $ x' $
en petit segments et on applique le raisonnement ci-dessus {\`a} chacun 
de ces petits segments.
\end{proof}
Nous utilisons les notations de \ref{SS::cycles}.\\
Notons $ G \subset H_{2(n-k)+1}(V \cap P)$ le sous groupe engendr{\'e} par les cycles $\g_1,\dots,\g_l$.
 Le lemme pr{\'e}c{\'e}dent et le th{\'e}or{\`e}me de section hyperplane
relatif (Proposition \ref{P::Lefschetz}) entra{\^\i}nent que
l'inclusion $V \cap P \subset V$ induit une application
surjective en homologie
$$\xymatrix{G \ar@{->>}[r]^-{\eta} &H_{2(n-k)+1}(V)} .$$
Autrement dit nous avons montr{\'e} que les cycles {\'e}vanescents engendre le groupe d'homologie $
H_{2(n-k)+1}(V) $.
Pour terminer la d{\'e}monstration, il nous reste {\`a}  montrer qu'ils
ne v{\'e}rifient pas de relation non triviale.\\
Soit $g=(g_1,\dots,g_k)$ un germe d'application pyramidale, consid{\'e}rons la propri{\'e}t{\'e} suivante:\\$(P_g)$ Pour tout $j \leq k$ le germe d'application $(g_1,\dots,g_j)$ est pyramidale.
\begin{lemme}{Il existe un ouvert de Zariski $\Omega$ de $GL(k,\CM)$ tel que pour
$R \in \Omega$, l'application $Rf$ v{\'e}rifie la propri{\'e}t{\'e} $P_{Rf}$.}
\end{lemme}
\begin{proof}
Une intersection finie d'ouvert de Zariski est un ouvert de Zariski;
par cons{\'e}quent il suffit de montrer la propri{\'e}t{\'e} pour $j=k-1$.\\
Notons $\t Z$ la sous vari{\'e}t{\'e} analytique de $GL(k,\CM) \times M$
qui consiste dans les paires $(R,x)$ pour lesquelles le germe de l'application
$R\bar f$ au point $x$ ne v{\'e}rifie pas la propri{\'e}t{\'e} $(P_{R\bar f})$.\\
Notons $C$ le lieu critique de $\bar f$ et $Z$ la projection de $Z$ dans
$GL(k,\CM)$.
On a le diagramme commutatif donn{\'e} par les projections sur le
premier et le deuxi{\`e}me facteur
$$\xymatrix{\t Z \subset GL(k,\CM) \times M\ar[r]^-{\pi_1}
                    \ar[d]^-{\pi_2}& M \supset C \\
                  Z \subset  GL(k,\CM)}$$
La condition de pyramidalit{\'e} entra{\^\i}ne que, au dessus d'un point
g{\'e}n{\'e}rique de $C$, la vari{\'e}t{\'e}  $\t Z$ est localement
isomorphe {\`a} $C \times Z$. Par ailleurs, la Proposition 3 entra{\^\i}ne
que au dessus d'un tel point $x \in M$,  $\dim \t Z < \dim C$ donc $Z$
est une sous vari{\'e}t{\'e} analytique de codimension au moins {\'e}gale {\`a} $1$.
Par cons{\'e}quent son compl{\'e}mentaire est un ouvert de Zariski. Le lemme
                    est d{\'e}montr{\'e}.
\end{proof}
Sans perte de g{\'e}n{\'e}ralit{\'e}, on peut donc supposer que l'application $f$ satisfait la
propri{\'e}t{\'e} $P_f$.\\
Notons $ A \subset M$
une fibre de Milnor de l'application pyramidale involutive $(\bar f_1,\dots,\bar f_{k-1})$, de telle sorte que $V
\subset A$ (si $k=1$, on prend $A=M$).\\
La diff{\'e}rence essentielle avec le cas des intersections compl{\`e}tes {\`a}
singularit{\'e}s isol{\'e}es provient du fait que $A$ est $2(n-k)+2$-connexe
(Th{\'e}or{\`e}me \ref{T::Lefschetz}).\\
Consid{\'e}rons la suite exacte du couple $ (A,V) $
$$\dots \to H_j(V) \to H_j(A) \to H_j(A,V) \to \dots .$$
Les {\'e}galit{\'e}s $ H_{2(n-k)+1}(A)=0,\ H_{2(n-k)+2}(A)=0 $ entra{\^\i}nent l'isomorphisme
$$
 H_{2(n-k)+1}(V) \simeq H_{2(n-k)+2}(A,V).$$
Soient $ D_1,\dots, D_l \subset (T \cap L)$ des voisinages ferm{\'e}s des valeurs
critiques de $\bar f$
et  $ \hat D_1,\dots,\hat D_l \subset A$ leurs pr{\'e}images par $\bar f $. Par
excision, on a un isomorphisme
$$ H_{2(n-k)+2}(A,V) \simeq \oplus_{i=1}^l H_{2(n-k)+2}(\hat D_i,V_i)$$

o{\`u} $ V_i $ est la fibre de $ f $ au dessus d'un point du bord de $
 D_i $.\\
Par excision, on a {\'e}galement l'isomorphisme
$$ H_{2(n-k)+2}(V,V \cap P) \simeq \oplus_{i=1}^l H_{2(n-k)+2}(\hat D_i \cap P,V_i \cap P)$$
o{\`u} $P$ est l'espace vectoriel qui a servi {\`a} construire les cycles
{\'e}vanescents (voir \ref{SS::cycles}).\\
Ces isomorphismes donnent lieu {\`a} un diagramme commutatif:
$$\xymatrix{ G \ \ \ar@{^{(}->}[r] \ar@{->>}@/^3.5pc/[rr]^{\eta} & H_{2(n-k)+1}(V \cap P)  \ar@{->>}[r]^-\xi& H_{2(n-k)+1}(V)  \\
 \ZM^l\  \ar@{^{(}->}[r]^-\nu  \ar@{->>}[u] \ar@/_2.8pc/[rr]^{\t \eta} & \oplus_{i=1}^l
 H_{2(n-k)+2}(\hat D_i \cap P,V_i \cap P) \ar[u] \ar[r] & \oplus_{i=1}^l H_{2(n-k)+2}(\hat D_i,V_i) \ar[u]_\simeq}$$
\vskip0.6cm
Les application $\nu$ et ${\t \eta}$ sont d{\'e}finies ainsi:
les cycles {\'e}vanescents qui engendrent le groupe $ G $
sont les classes d'homologie des bords boules ferm{\'e}es
$ C_i \subset \hat D_i \cap P $,\
$\d C_i \simeq S^{2(n-k)+1}$, $i=1,\dots,l$; les applications $\nu$
 en ${\t \eta}$ envoient
 un {\'e}l{\'e}ment du groupe ab{\'e}lien libre engendr{\'e} par les $C_i$ sur sa
 classe d'homologie respective.\\
L'application $\eta$ est surjective, par cons{\'e}quent
${\t \eta}$ l'est {\'e}galement.\\
En prenant pour syst{\`e}me g{\'e}n{\'e}rateur de l'image de ${\t \eta}$ les classes
 d'homologie des $C_i$, l'application ${\t \eta} $ est diagonale.
Il reste donc {\`a} montrer que pour tout $i \leq l$ et pour
tout $n \in \ZM$, l'image par ${\t \eta}$ de la classe de $nC_i$ est non nulle.\\
Ce dernier point provient du fait que
 le volume de  la boule $ C_i $ est non nul. En effet, si l'on note $ \a=\sum_{i=1}^n p_idq_i $ la forme d'action, on a par
Stokes
$$\int_{\d C_i}\a \w (d\a)^{n-k}=\int_{C_i}(d\a)^{(n-k)+1} \neq 0  $$
donc la classe de $n \d C_i $ dans $ H_{2(n-k)+1}(V) $ est non nulle. Il 
d{\'e}coule par cons{\'e}quent de l'isomorphisme $ H_{2(n-k)+1}(V) \simeq 
H_{2(n-k)+2}(A,V) $ que la classe de $n C_i $ dans $ H_{2(n-k)+2}(\hat
D_i,V_i) \subset H_{2(n-k)+2}(A,V) $  est non nulle. Ce qui prouve
que l'application ${\t \eta}$ est un isomorphisme et ach{\`e}ve la d{\'e}monstration du th{\'e}or{\`e}me.\\
 
{\em Remerciements.}{ Je remercie D.T. L{\^e},
F. Pham, B. Teissier et D. van Straten
pour les explications qu'ils m'ont donn{\'e}es sur divers points de la
th{\'e}orie des singularit{\'e}s ainsi que V.I. Arnold qui m'a enseign{\'e} la th{\'e}orie de
Picard-Lefschetz alors que j'{\'e}tais {\'e}tudiant. Je remercie le r{\'e}f{\'e}r{\'e}e
de m'avoir signal{\'e} plusieurs impr{\'e}cisions dans la version initiale
de travail.\\ Enfin, je remercie le Deutsche Forschungsgemeinschaft pour avoir financ{\'e}
ce travail ainsi que l'Institut des Hautes
{\'E}tudes Scientifiques pour m'avoir accueilli d'Octobre {\`a} D{\'e}cembre 2004
dans le cadre du 6{\`e}me programme Europ{\'e}en de la Commission Europ{\'e}enne (Contrat RITA-CT-2004-505493).}
\bibliographystyle{amsplain}
\bibliography{master}
\end{document}